\documentclass[graybox]{svmult}

\usepackage{helvet}         
\usepackage{courier}        
\usepackage{type1cm}        
\usepackage{makeidx}         
\usepackage{graphicx}        
\usepackage{multicol}        
\usepackage[bottom]{footmisc}

\makeindex             

\usepackage{amsfonts}
\usepackage{epsfig}
\usepackage{latexsym, amssymb, color, calc, amsmath, dsfont}
\usepackage{mathrsfs}
\usepackage{verbatim}
\usepackage[caption=false]{subfig}
\usepackage{common_symbols}
\usepackage{afterpage}
\graphicspath{{img/}}

\newcommand{\reviewerA}[1]{{#1}}

\renewcommand{\u}{\ensuremath{\vector{u}}}
\renewcommand{\v}{\ensuremath{\vector{v}}}
\newcommand{\z}{\ensuremath{\vector{z}}}
\newcommand{\p}{\ensuremath{p}}

\newcommand{\n}{\ensuremath{\vector{n}}}
\newcommand{\x}{\ensuremath{\vector{x}}}

\begin{document}

\title*{Reduced-order semi-implicit schemes for fluid-structure interaction problems}
\titlerunning{Reduced-order semi-implicit schemes for FSI}
\author{Francesco Ballarin, Gianluigi Rozza and Yvon Maday}
\institute{Francesco Ballarin \at Mathematics Area, mathLab, SISSA, via Bonomea 265, I-34136 Trieste, Italy, \email{francesco.ballarin@sissa.it}
\and Gianluigi Rozza \at Mathematics Area, mathLab, SISSA, via Bonomea 265, I-34136 Trieste, Italy, \email{gianluigi.rozza@sissa.it}
\and Yvon Maday \at Sorbonne Universit\'{e}s, UPMC Universit\'{e} Paris 06 and CNRS UMR 7598, Laboratoire Jacques-Louis Lions, F-75005, Paris, France, \email{maday@ann.jussieu.fr} \at Division of Applied Mathematics, Brown University, Providence RI, USA
}
\maketitle

\abstract{POD--Galerkin reduced-order models (ROMs) for fluid-structure interaction problems (incompressible fluid and thin structure) \reviewerA{are} proposed in this paper. Both the high-fidelity and reduced-order methods are based on a Chorin-Temam operator-splitting approach. Two different reduced-order methods are proposed, which differ on velocity continuity condition, imposed weakly or strongly, respectively. The resulting ROMs are tested and compared on a representative haemodynamics test case characterized by wave propagation, in order to assess the capabilities of the proposed strategies.}

\section{Introduction}
Several applications are characterized by multi-physics phenomena, such as the interaction between an incompressible fluid and a compressible structure. The capability to perform real-time multi-physics simulations could greatly increase the applicability of computational methods in applied sciences and engineering. To reach this goal, reduced-order modelling techniques are applied in this paper. We refer the interested reader to \cite{Amsallem2010,LassilaQuarteroniRozza2012,Colciago2014,BallarinRozza2016} for some representative previous approaches to the reduction of fluid-structure interaction problems, arising in aeroelasticity \cite{Amsallem2010} or haemodynamics \cite{LassilaQuarteroniRozza2012,Colciago2014,BallarinRozza2016}. The reduction proposed in the current work is based on a POD--Galerkin approach. A difference with our previous work \cite{BallarinRozza2016}, where a reduced-order monolithic approach has been proposed, is related to the use of a partitioned reduced-order model, based on the semi-implicit operator-spitting approach originally employed in \cite{Fernandez2007} for the high-fidelity method. Extension to other methods (e.g. \cite{BadiaQuaini2008,Guidoboni2009}) is possible and object of forthcoming work. The formulation of the FSI problem is summarized in Section \ref{sec:problem}, and its high-fidelity discretization is reported in Section \ref{sec:HF}. Two reduced-order models are proposed in Section \ref{sec:ROM}, and compared by means of a numerical test case in Section \ref{sec:numerical_comparison}. Conclusions and perspectives are summarized in the final section of the paper.

\section{Problem formulation}
\label{sec:problem}
In this section the formulation of the fluid-structure interaction (FSI) model problem is summarized. Let us consider the bidimensional fluid domain $\Omega = [0, L] \times [0, h_f]$. Its boundary is composed of a compliant wall $\Sigma = [0, L] \times \{h_f\}$ (top), fluid inlet section $\Gamma_{in} = \{0\} \times [0, h_f]$ (left) and fluid outlet section $\Gamma_{out} = \{L\} \times [0, h_f]$ (right), and a wall $\Gamma_{sym} = [0, L] \times \{0\}$ (bottom). For the sake of simplicity we assume low Reynolds numbers for the fluid problem and infinitesimal displacements for the compliant wall. Thus, in the following, we will consider \emph{unsteady Stokes} equations on a \emph{fixed} domain $\Omega$ for the fluid and a linear structural model for the compliant wall. In particular, we will further assume that the structure undergoes negligible horizontal displacements, so that the structural equations of the compliant wall can be described by a \emph{generalized string model} \cite{QuarteroniFormaggia2004,QuarteroniTuveriVeneziani2000}.

The coupled fluid-structure interaction problem is therefore: for all $t \in (0, T]$, find fluid velocity $\u(t): \Omega \to \RR^2$, fluid pressure $p(t): \Omega \to \RR$ and structure displacement $\eta(t): \Sigma \to \RR$ such that
\begin{equation}
\begin{cases}
\rho_f \partial_t \u - \div(\boldsymbol{\sigma}(\u, p)) = \vector{0} & \text{ in } \Omega \times (0, T],\\
\div \u = 0 &  \text{ in } \Omega  \times (0, T],\\
\u = \partial_t \eta\, \n & \text{ on } \Sigma  \times (0, T],\\
\rho_s h_s \partial_{tt} \eta - c_1 \partial_{xx} \eta + c_0 \eta = - \boldsymbol{\sigma}(\u, p) \n \cdot \n & \text{ on } \Sigma  \times (0, T].
\end{cases}
\label{eq:FSI}
\end{equation}
Equations \eqref{eq:FSI}$_1$-\eqref{eq:FSI}$_2$ formulate the Stokes problem on the fixed fluid domain $\Omega$, having defined the fluid Cauchy stress tensor $\boldsymbol{\sigma}(\u, p) := - p \boldsymbol{I} + 2 \mu_f \boldsymbol{\varepsilon}(\u)$, $\boldsymbol{\varepsilon}(\u) := \frac{1}{2} (\grad\u + \grad^T \u)$, while \eqref{eq:FSI}$_4$ is the equation of the structural motion of the compliant wall. Continuity of the velocity on the interface is guaranteed by \eqref{eq:FSI}$_3$. The FSI system is completed by suitable initial and boundary conditions. In this paper we assume resting conditions at $t = 0$ and the following fluid boundary conditions on $\partial\Omega \setminus \Sigma$:
\begin{equation}
\begin{cases}
\boldsymbol{\sigma}(\u, p) \n = - p_{in}(t) \n & \text{ on } \Gamma_{in}  \times (0, T],\\
\boldsymbol{\sigma}(\u, p) \n = - p_{out}(t) \n & \text{ on } \Gamma_{out}  \times (0, T],\\
\u \cdot \n = 0, \quad \boldsymbol{\sigma}(\u, p) \n \cdot \vector{\tau} = 0 & \text{ on } \Gamma_{sym}  \times (0, T],\\
\end{cases}
\label{eq:bcFluid}
\end{equation}
and the following structure boundary condition on $\partial\Sigma$:
\begin{equation}
\eta = 0 \quad \text{ on } \partial\Sigma  \times (0, T].
\label{eq:bcStruct}
\end{equation}
Equations \eqref{eq:bcFluid}$_1$ and \eqref{eq:bcFluid}$_2$ prescribe pressure on inlet and outlet section, respectively; equation \eqref{eq:bcFluid}$_3$ is a symmetry condition, arising from the consideration of this simplified 2D problem as a section of a 3D cylindrical configuration. Here $\n$ and $\vector{\tau}$ denote the outer unit normal to $\Omega$ and tangential vector to $\Sigma$, respectively. Finally, \eqref{eq:bcStruct} prescribes a clamped wall near both the inlet and outlet section of the fluid. The value of constitutive parameters $\rho_f$ (fluid density), $\rho_s$ (structure density), $\mu_f$ (fluid viscosity), $c_1$ and $c_0$ (structure constitutive parameters), $L$ (domain width), $h_f$ (fluid height), $h_s$ (structure thickness) will be further specified in Section \ref{sec:numerical_comparison}, as well as choices of $p_{in}(t)$, $p_{out}(t)$ and $T$ for which the FSI system \eqref{eq:FSI}-\eqref{eq:bcStruct} is characterized by propagation of a pressure wave. 

\section{High-fidelity formulation: semi-implicit scheme}
\label{sec:HF}
In this section we summarize the high-fidelity discretization of the FSI system \eqref{eq:FSI}-\eqref{eq:bcStruct}. An operator splitting approach, based on a Chorin-Temam projection scheme, is pursued. In particular, Robin-Neumann iterations are carried out in order to enhance the stability of the resulting algorithm.

\subsection{A projection-based semi-implicit coupling scheme}
We employ in this work a projection-based semi-implicit scheme, as proposed in \cite{Fernandez2007,Astorino2010}. The fluid equations are discretized in time using the Chorin-Temam projection scheme \cite{QuarteroniValli2008}. Thus, denoting by $\Delta t$ the time step length, $D_t f^{k+1} := \frac{f^{k+1} - f^{k}}{\Delta t}$ the first backward difference approximation of the time derivative of $f(t^{n+1})$ and $D_{tt} f^{k+1} = D_t(D_t f^{k+1})$, we consider the following semi-implicit time discretization of \eqref{eq:FSI}-\eqref{eq:bcStruct}: for any $k = 1, \hdots, K = T/\Delta t$
\begin{enumerate}
\item Explicit step (fluid viscous part): find\footnote{For the sake of an easier comparison with Remark 1 we employ here the $\widetilde{\cdot}$ notation for the velocity. However, since we will always employ the pressure Poisson formulation, the $\widetilde{\cdot}$ will be dropped in the next sections.} $\widetilde{\u}^{k+1}: \Omega \to \RR^2$ such that:
\begin{equation}
\begin{cases}
\rho_f \frac{\widetilde{\u}^{k+1} - \widetilde{\u}^{k}}{\Delta t} - 2 \mu_f \div\, \boldsymbol{\varepsilon}(\widetilde{\u}^{k+1}) = - \grad p^k & \text{ in } \Omega,\\
\widetilde{\u}^{k+1} = D_t \eta^k\, \n & \text{ on } \Sigma.
\end{cases}
\label{eq:si_1}
\end{equation}
\item Implicit step:
\begin{enumerate}
\item[2.1.] Fluid projection substep: find $p^{k+1}: \Omega \to \RR$ such that:
\begin{equation}
\begin{cases}
- \div(\grad p^{k+1}) = - \frac{\rho_f}{\Delta t} \div \widetilde{\u}^{k+1} & \text{ in } \Omega,\\
\frac{\partial}{\partial \n} \p^{k+1} = - \rho_f D_{tt} \eta^{k+1} & \text{ on } \Sigma.
\end{cases}
\label{eq:si_21}
\end{equation}
\item[2.2.] Structure substep: find $\eta^{k+1}: \Sigma \to \RR$ such that:
\begin{equation}
\rho_s h_s D_{tt} \eta^{k+1} - c_1 \partial_{xx} \eta^{k+1} + c_0 \eta^{k+1} = - \boldsymbol{\sigma}(\widetilde{\u}^{k+1}, p^{k+1}) \n \cdot \n \quad \text{ on } \Sigma.
\label{eq:si_22}
\end{equation}
\end{enumerate}
\end{enumerate}

The implicit step couples pressure stresses to the structure, and it is iterated until convergence.

\begin{remark}
\label{rem:pressureDarcy}
In place of step 1 and 2.1 (pressure Poisson formulation) one could also consider the following pressure Darcy formulation:
\begin{enumerate}
\item[I.] Explicit step (fluid viscous part): find $\widetilde{\u}^{k+1}: \Omega \to \RR^2$ such that:
\begin{equation*}
\begin{cases}
\rho_f \frac{\widetilde{\u}^{k+1} - \u^{k}}{\Delta t} - 2 \mu_f \div\, \boldsymbol{\varepsilon}(\widetilde{\u}^{k+1}) = \vector{0} & \text{ in } \Omega,\\
\widetilde{\u}^{k+1} = D_t \eta^k\, \n & \text{ on } \Sigma.
\end{cases}
\end{equation*}
\item[II.] Implicit step:
\begin{enumerate}
\item[II.1.] Fluid projection substep: find $\u^{k+1}: \Omega \to \RR^2$ and $p^{k+1}: \Omega \to \RR$ such that:
\begin{equation*}
\begin{cases}
\rho_f \frac{\u^{k+1} - \widetilde{\u}^{k+1}}{\Delta t} + \grad p^{k+1} = \vector{0} & \text{ in } \Omega,\\
\u^{k+1} \cdot \n = D_t \eta^{k+1} & \text{ on } \Sigma.
\end{cases}
\end{equation*}
\item[II.2.] Structure substep: find $\eta^{k+1}: \Sigma \to \RR$ such that:
\begin{equation*}
\rho_s h_s D_{tt} \eta^{k+1} - c_1 \partial_{xx} \eta^{k+1} + c_0 \eta^{k+1} = - \boldsymbol{\sigma}(\widetilde{\u}^{k+1}, p^{k+1}) \n \cdot \n \quad \text{ on } \Sigma.
\end{equation*}
\end{enumerate}
\end{enumerate}

For the sake of a more efficient reduced-order model (see Section \ref{sec:ROM}) it is convenient to consider the pressure Poisson formulation (steps 1 and 2.1) rather than the pressure Darcy formulation (step I and II.1), \reviewerA{because the latter would require a larger system (comprised of both velocity and pressure) at step II.1}.
\end{remark}

Finally, in order to enhance the stability of the projection method we employ a Robin-Neumann coupling, as proposed in \cite{Astorino2010}. See also \cite{BadiaVergara2008,Fernandez2016} for related topics. Thus, we replace \eqref{eq:si_21}$_2$ with 
\begin{equation}
\frac{\partial}{\partial \n} \p^{k+1} + \alpha_{Rob} \p^{k+1} = - \rho_f D_{tt} \eta^{k+1} + \alpha_{Rob} \p^{k, *} \quad \text{ on } \Sigma.
\end{equation}
being $\alpha_{Rob} > 0$ and $\p^{k, *}$ an extrapolation of the pressure \reviewerA{(which will be defined in section 3.2)}. 
In particular, following \cite{Fernandez2013a,Fernandez2013b}, we choose $\alpha_{Rob} := \frac{\rho_f}{\rho_s h_s}$.
\reviewerA{We remark that, due to the simplifying assumptions of this problem (linear structural model, fixed domain), the implicit step could have been solved in one shot, since it defines a linear system in $(p^{k+1}, \eta^{k+1})$. We still keep the Robin-Neumann coupling in this case in order to assess the capabilities of such procedure in a reduced-order setting, since it will be required in more general nonlinear problems.}

\subsection{Space discretization of the high-fidelity formulation}
Denote by $V = [H^1(\Omega)]^2$ the fluid velocity space (endowed with the $H^1$ seminorm), by $Q=L^2(\Omega)$ the fluid pressure space (endowed with the $L^2$ norm), and by $E=H^1(\Sigma)$ the structure displacement space (endowed with the $H^1$ seminorm). After having obtained a weak formulation of the semi-implicit formulation, we consider a finite element (FE) discretization for steps 1, 2.1 and 2.2. Second order Lagrange FE are employed for fluid velocity (step 1) and structural displacement (step 2.2), resulting in FE spaces $V_h \subset V$ and $E_h \subset E$, respectively, while fluid pressure is discretized by first order Lagrange FE, $Q_h \subset Q$. Thus, the corresponding Galerkin-FE formulation reads: for any $k = 1, \hdots, K$
\begin{enumerate}
\item[1$_h$.] Explicit step (fluid viscous part): find ${\u}^{k+1}_h \in V_h$ such that:
\begin{equation}
\int_\Omega \frac{\rho_f}{\Delta t} \u^{k+1}_h \cdot \v_h \,d\x \,+ \int_\Omega 2 \mu_f \boldsymbol{\varepsilon}({\u}^{k+1}_h) : \grad\v_h \,d\x = 
\int_\Omega \frac{\rho_f}{\Delta t} \u^{k}_h \cdot \v_h \,d\x \, - \int_\Omega \grad p^k_h \cdot \v_h \,d\x
\label{eq:si_1_h}
\end{equation}
for all $\v_h \in V_h$, subject to the coupling condition
\begin{equation}
{\u}^{k+1}_h = D_t \eta^k_h\, \n \quad \text{ on } \Sigma  \times (0, T],
\label{eq:bc_strongly}
\end{equation}
and to the boundary condition ${\u}^{k+1}_h \cdot \n = 0$ on $\Gamma_{sym}  \times (0, T]$.
\item[2$_h$.] Implicit step: for any $j = 0, \hdots$, until convergence:
\begin{enumerate}
\item[2.1$_h$.] Fluid projection substep, with Robin boundary conditions: find $p^{k+1,j+1}_h \in Q_h$ such that:
\begin{align*}
\int_\Omega \grad p^{k+1,j+1}_h \cdot \grad q_h \,d\x &+ \int_\Sigma \alpha_{Rob}\, p^{k+1,j+1}_h \, q_h \,ds = \\
& - \int_\Omega \frac{\rho_f}{\Delta t} \div {\u}^{k+1}_h \, q_h  \,d\x
 - \int_\Sigma \rho_f D_{tt} \eta^{k+1,j}_h \, q_h \,ds \\
& + \int_\Sigma \alpha_{Rob}\, p^{k+1,j}_h \, q_h \,ds
\end{align*}
for all $q_h \in Q_h$, subject to the boundary conditions $p^{k+1,j+1}_h = p_{in}(t)$ on $\Gamma_{in} \times (0, T]$ and $p^{k+1,j+1}_h = p_{out}(t)$ on $\Gamma_{out} \times (0, T]$. \reviewerA{Here the value $p^{k+1,j}_h$ has been chosen as pressure extrapolation for the Robin-Neumann coupling.}
\item[2.2$_h$.] Structure substep: find $\eta^{k+1,j+1}_h \in E_h$ such that:
\begin{align}
\int_\Sigma \frac{\rho_s h_s}{\Delta t^2} \, \eta^{k+1,j+1}_h \, \zeta_h \,ds
& + \int_\Sigma c_1 \partial_{x} \eta^{k+1,j+1}_h \, \partial_x \zeta_h \,ds \notag\\
&+ \int_\Sigma c_0 \eta^{k+1,j+1}_h \, \zeta_h \,ds = 
\int_\Sigma \frac{\rho_s h_s}{\Delta t^2} \eta^{k}_h \, \zeta_h \,ds \label{eq:si_22_h}\\
& + \int_\Sigma \frac{\rho_s h_s}{\Delta t} D_t \eta^{k}_h \, \zeta_h \,ds
- \int_\Sigma \boldsymbol{\sigma}({\u}^{k+1}, p^{k+1,j+1}) \n \cdot \zeta_h \n \,ds\notag
\end{align}
for all $\zeta_h \in E_h$, subject to the boundary conditions $\eta^{k+1,j+1}_h = 0$ on $\partial\Sigma$.
\end{enumerate}
\end{enumerate}

The coupling condition \eqref{eq:bc_strongly} is imposed strongly. We will further comment  in Section \ref{sec:ROM} on the imposition of this condition at the reduced-order level. A relative error on the increments is chosen as stopping criterion for step $2_h$, that is the implicit step is repeated until
\begin{equation*}
\min\left\{\frac{\Norm{p^{k+1,j+1}_h - p^{k+1,j}_h}_{Q}}{\Norm{p^{k+1,j+1}_h}_{Q}}, \frac{\Norm{\eta^{k+1,j+1}_h - \eta^{k+1,j}_h}_{E}}{\Norm{\eta^{k+1,j+1}_h}_{E}}\right\} < \text{tol},
\end{equation*}
for some prescribed tolerance $\text{tol}$. The solution $(p^{k+1,j^*}_h, \eta^{k+1,j^*}_h)$ at the iteration $j^*$ such that convergence is achieved is then denoted by $(p^{k+1}_h, \eta^{k+1}_h)$.

\section{Reduced-order formulation: POD--Galerkin semi-implicit scheme}
\label{sec:ROM}

In this section we propose two Proper Orthogonal Decomposition (POD)--Galerkin semi-implicit reduced-order models (ROMs) for FSI system \eqref{eq:FSI}-\eqref{eq:bcStruct}. The first ROM (FSI ROM 1) is built starting directly from steps $1_h$, $2.1_h$ and $2.2_h$, and performing a Galerkin projection. Special treatment will be devoted to the imposition of the coupling condition \eqref{eq:bc_strongly}; unfortunately, this requires enlarging the reduced-order systems. The second ROM (FSI ROM 2) will exploit a simple change of variable for the fluid velocity to bypass this issue. 
In both cases, an offline-online computational decoupling is sought \cite{RozzaHuynhPatera2007}.

\subsection{FSI ROM 1 approach}
\subsubsection{Offline stage}
\label{sec:ROM_1_offline}
During the offline stage, the solution of the high-fidelity problem $1_h$, $2.1_h$ and $2.2_h$ is computed. We then consider the following snapshot matrices
\begin{align*}
S_{\u} = [\underline{\mathbf{u}}^{1}_h | \hdots | \underline{\mathbf{u}}^{K}_h] \in \RR^{N_h^{\u} \times K},\\
S_{\p} = [\underline{\mathbf{p}}^{1}_h | \hdots | \underline{\mathbf{p}}^{K}_h] \in \RR^{N_h^{\p} \times K},\\
S_{\eta} = [\underline{\boldsymbol{\eta}}^{1}_h | \hdots | \underline{\boldsymbol{\eta}}^{K}_h] \in \RR^{N_h^{\eta} \times K},
\end{align*}
where we denote with the underlined notation the vector of FE degrees of freedom corresponding to each solution. Here $N_h^{\u} = \text{dim}(V_h)$, $N_h^{\p} = \text{dim}(Q_h)$ and $N_h^{\eta} = \text{dim}(E_h)$. Then, we carry out a \reviewerA{proper orthogonal decomposition} of each snapshot matrix; \reviewerA{the method of snapshots is used, and the snapshots are weighted with the inner product associated to their functional space.} Then, the first $N^{\u}$, $N^{\p}$ and $N^{\eta}$ (respectively) left singular vectors, denoted by $\{\vector{\varphi}_i\}_{i = 1}^{N^{\u}}$, $\{\psi_j\}_{j = 1}^{N^{\p}}$ and $\{\phi_l\}_{l = 1}^{N^{\eta}}$ (resp.), are then chosen as basis functions for the reduced spaces $V_N^{(1)}$, $Q_N^{(1)}$ and $E_N^{(1)}$ (resp.), i.e.\
\begin{equation*}
V_N^{(1)} = \text{span}(\{\vector{\varphi}_i\}_{i = 1}^{N^{\u}}), \quad
Q_N^{(1)} = \text{span}(\{\psi_j\}_{j = 1}^{N^{\p}}), \quad
E_N^{(1)} = \text{span}(\{\phi_l\}_{l = 1}^{N^{\eta}}).
\end{equation*}

\subsubsection{On the imposition of coupling condition \eqref{eq:bc_strongly}}
\label{sec:ROM_1_mult}
The major drawback of this approach is related to the fact that the reduced spaces $V_N^{(1)}$ and $E_N^{(1)}$ do not guarantee, in general, that the coupling conditions \eqref{eq:bc_strongly} holds. To this end, during the online stage, we will resort to a weak imposition of \eqref{eq:bc_strongly} by Lagrange multipliers. More precisely, we will enforce weakly the normal component of \eqref{eq:bc_strongly} (i.e. ${\u}^{k+1}_h \cdot \n = D_t \eta^k_h$), while the tangential component of \eqref{eq:bc_strongly} (i.e. ${\u}^{k+1}_h \cdot \vector{\tau} = 0$) is already imposed strongly, since it is homogeneous and all basis functions in $V_N$ satisfy it. Thus, during the offline stage we need to build an additional snapshot matrix of the Lagrange multipliers, in order to carry out a POD to obtain a reduced Lagrange multipliers space. The traction on the interface, which can be evaluated as the residual of \eqref{eq:si_1_h} for test functions $\v_h := v_h \n$ which do not vanish on the interface, is indeed the Lagrange multiplier to \eqref{eq:bc_strongly}. Therefore, we build an additional snapshot matrix
\begin{align*}
S_{\lambda} = [\underline{\boldsymbol{\lambda}}^{1}_h | \hdots | \underline{\boldsymbol{\lambda}}^{K}_h] \in \RR^{N_h^{\u} \times K},
\end{align*}
containing the FE degrees of freedom corresponding to the residual of \eqref{eq:si_1_h} for test functions $\v_h := v_h \n$ that do not vanish on the interface, compute a POD, and, similarly to the previous section, obtain a reduced space $L_N^{(1)}$ as the space spanned by the first $N_{\vector{\lambda}}$ left singular vectors. \reviewerA{During the POD computation by the method of snapshots we employ the $L^2$ inner product on the interface as weight.} We remark again here that the Lagrange multiplier approach is \emph{not} actually used during the high-fidelity solution of the FSI system (in favor of a strong imposition), but rather the snapshot matrix $S_{\lambda}$ is obtained as a post-processing of the obtained solution. In contrast, in the online stage the Lagrange multiplier approach will actually be used while solving the linear system associated to the reduced fluid viscous step, in order to impose weakly the coupling condition \eqref{eq:bc_strongly}.

\subsubsection{Online stage}
\label{sec:ROM_1_online}
A reduced-order approximation of the FSI problem is then obtained by means of a Galerkin projection over the reduced spaces $V_N^{(1)}, Q_N^{(1)}$ and $E_N^{(1)}$, respectively, treating the coupling condition \eqref{eq:bc_strongly} with Lagrange multipliers in the reduced space $L_N^{(1)}$. Thus, the corresponding online stage of the POD--Galerkin method reads: for any $k = 1, \hdots, K$
\begin{enumerate}
\item[1$_N^{(1)}$.] Explicit step (fluid viscous part), with weak imposition of coupling conditions through Lagrange multipliers: find $({\u}^{k+1}_N, \lambda^{k+1}_N) \in V_N^{(1)} \times L_N^{(1)}$ such that:
\begin{align*}
\begin{cases}
\int_\Omega \frac{\rho_f}{\Delta t} \u^{k+1}_N \cdot \v_N \,d\x & + \int_\Omega 2 \mu_f \boldsymbol{\varepsilon}({\u}^{k+1}_N) : \grad\v_N \,d\x \\
& + \int_\Sigma \lambda^{k+1}_N \n \cdot \v_N \,ds = 
\int_\Omega \frac{\rho_f}{\Delta t} \u^{k}_N \cdot \v_N \,d\x \, \\
& - \int_\Omega \grad p^k_N \cdot \v_N \,d\x,\\
\int_\Sigma {\u}^{k+1}_N \cdot \Upsilon_N \n \,ds & = \int_\Sigma D_t \eta^k_h\,\, \Upsilon_N \,ds,
\end{cases}
\end{align*}
for all $(\v_N, \Upsilon_N) \in V_N^{(1)} \times L_N^{(1)}$. We note that the boundary condition ${\u}^{k+1}_N \cdot \n = 0$ on $\Gamma_{sym}$ is implicitly verified, since it is satisfied by any element in $V_N$.
\item[2$_N^{(1)}$.] Implicit step: for any $j = 0, \hdots$, until convergence:
\begin{enumerate}
\item[2.1$_N^{(1)}$.] Fluid projection substep, with Robin boundary conditions: find $p^{k+1,j+1}_N \in Q_N^{(1)}$ such that:
\begin{align*}
\int_\Omega \grad p^{k+1,j+1}_N \cdot \grad q_N \,d\x &+ \int_\Sigma \alpha_{Rob}\, p^{k+1,j+1}_N \, q_N \,ds = \\
& - \int_\Omega \frac{\rho_f}{\Delta t} \div {\u}^{k+1}_N \, q_N  \,d\x
 - \int_\Sigma \rho_f D_{tt} \eta^{k+1,j}_N \, q_N \,ds \\
& + \int_\Sigma \alpha_{Rob}\, p^{k+1,j}_N \, q_N \,ds
\end{align*}
for all $q_N \in Q_N^{(1)}$. The imposition of the boundary conditions $p^{k+1,j+1}_N = p_{in}(t)$ on $\Gamma_{in}$ and $p^{k+1,j+1}_N = p_{out}(t)$ on $\Gamma_{out}$, although not automatically prescribed by the reduced space $Q_N^{(1)}$, can be easily treated by a lifting in $2.1_h$ without the need to introduce an additional Lagrange multiplier for the pressure, since the values to be imposed value do not depend on any reduced space, rather are known functions. The details are omitted for the sake of brevity; the interested reader is referred to \cite{BallarinManzoniQuarteroniRozza2015} for more details.
\item[2.2$_N^{(1)}$.] Structure substep: find $\eta^{k+1,j+1}_N \in E_N^{(1)}$ such that:
\begin{align*}
\int_\Sigma \frac{\rho_s h_s}{\Delta t^2} \, \eta^{k+1,j+1}_N \, \zeta_N \,ds
& + \int_\Sigma c_1 \partial_{x} \eta^{k+1,j+1}_N \, \partial_x \zeta_N \,ds \\
&+ \int_\Sigma c_0 \eta^{k+1,j+1}_N \, \zeta_N \,ds = 
\int_\Sigma \frac{\rho_s h_s}{\Delta t^2} \eta^{k}_N \, \zeta_N \,ds \\
& + \int_\Sigma \frac{\rho_s h_s}{\Delta t} D_t \eta^{k}_N \, \zeta_N \,ds
- \int_\Sigma \boldsymbol{\sigma}({\u}^{k+1}, p^{k+1,j+1}) \n \cdot \zeta_N \n \,ds
\end{align*}
for all $\zeta_N \in E_N^{(1)}$. The boundary condition $\eta^{k+1,j+1}_N = 0$ on $\partial\Sigma$ is implicitly verified.
\end{enumerate}
\end{enumerate}

As for the high-fidelity model, a stopping criterion on the relative increment of the solution is employed to terminate step $2_N^{(1)}$.

\begin{remark}[On efficient offline-online decoupling]
\label{rem:offline_online_decoupling_1}
The reduced-order problem $1_N^{(1)}$, $2.1_N^{(1)}$ and $2.2_N^{(1)}$ can easily account for an efficient offline-online decoupling, thanks to the linearity assumption in the problem formulation. For instance, the fluid mass term $\int_\Omega \frac{\rho_f}{\Delta t} \u^{k+1}_N \cdot \v_N \,d\x$ in $1_N^{(1)}$, is efficiently assembled at the end of the offline stage as
\begin{equation*}
M_N^{(1)} := (Z_N^{\u})^T \, M_h \, Z_N^{\u},
\end{equation*}
and loaded during the online stage. Here $Z_N^{\u}$ is the matrix which contains the velocity basis functions $\{\vector{\varphi}_i\}_{i = 1}^{N^{\u}}$ as columns, and $M_h$ is the FE matrix corresponding to fluid mass term $\int_\Omega \frac{\rho_f}{\Delta t} \u^{k+1}_h \cdot \v_h \,d\x$ in $1_h$. One can carry out a similar computational procedure for all terms in the reduced formulation $1_N^{(1)}$, $2.1_N^{(1)}$ and $2.2_N^{(1)}$.

In a more general (nonlinear, geometrical parametrized) setting one can resort to the empirical interpolation method \cite{BarraultMadayNguyenPatera2004} to recover an efficient offline-online splitting, as recently shown for FSI problems in \cite{BallarinRozza2016}.
\end{remark}

\begin{remark}[On supremizer enrichment: \reviewerA{the role of pressure}]
In contrast to what is usually done in the reduced basis approximation of parametrized fluid dynamics problem (see \cite{RozzaVeroy2007,Rozza2013,Rozza2009}, and also \cite{BallarinManzoniQuarteroniRozza2015} for an extension to POD--Galerkin methods) and previous works on FSI \cite{BallarinRozza2016,Colciago2014,LassilaManzoniQuarteroniRozza2013a}, we do not employ in this case a supremizer enrichment of the velocity space \reviewerA{to enforce inf-sup stability of the mixed velocity-pressure formulation}. This is heuristically motivated by the fact that, even at the high-fidelity level, the Chorin-Temam scheme, in its pressure Poisson version, can be successfully applied to FE spaces that do not fulfill a $(V, Q)$-inf-sup condition \cite{guermond1998stability}, even though it may result in non-optimal error estimates. 
\end{remark}

\begin{remark}[On supremizer enrichment: \reviewerA{the role of Lagrange multiplier}]
\reviewerA{%
Problem $1_N^{(1)}$ still features a saddle point structure. We remark that this structure is not due to the original problem, but rather due to our choice of coupling conditions by Lagrange multipliers. A drawback of ROM 1 is now apparent for what concerns the size of the reduced system $1_N^{(1)}$, which needs to be increased to $N_{\u} + N_{\vector{\lambda}}$ due to weak imposition of coupling condition. The ROM proposed in the next section has been devised to overcome this limitation, and results in a reduced explicit step of size $N_{\u}$.
Moreover, a further increase in dimension would be required if we were willing to enrich the velocity space $V_N^{(1)}$ with supremizers corresponding to  the inf-sup condition associated to problem $1_N^{(1)}$, i.e. solutions to
\begin{equation*}
\int_\Omega \nabla\vector{s}^k \cdot \nabla\vector{v} = \int_{\Sigma} \vector{\lambda}^k_h \cdot \vector{v}  \quad \forall \vector{v} \in V,\\
\end{equation*}
for all $k = 1, \hdots, K$. In this work we do not carry out such enrichment since it would further increase the size of the reduced explicit step; nevertheless, a detailed investigation of the stability of $1_N^{(1)}$ with and without enrichment by $\vector{s}^k$ is an ongoing task and will be presented in a forthcoming work.%
}
\end{remark}

\subsection{FSI ROM 2 approach}
\label{sec:ROM_2}
As we have seen in the previous section, it is challenging to enforce \eqref{eq:bc_strongly} at the reduced-order level. The second reduced-order model proposed in this paper overcomes these difficulties performing a change of variable for the fluid velocity, namely defining an auxiliary unknown $\vector{z}^{k+1}: \Omega \to \RR^2$ as
\begin{equation}
\vector{z}^{k+1} = \u^{k+1} - D_t \widehat{\eta}^k\, \n,
\label{eq:def_z}
\end{equation}
where $\widehat{\eta}^k$ is the solution of the following \emph{harmonic extension} problem
\begin{equation*}
- \Delta \widehat{\eta}^k = 0 \quad \text{ in } \Omega,
\end{equation*}
subject to the following inhomogeneous boundary condition on the interface
\begin{equation*}
\widehat{\eta}^k = \eta^k \quad \text{ on } \Sigma,
\end{equation*}
and homogeneous boundary condition on the remaining boundaries. In this way, the coupling condition \eqref{eq:bc_strongly} is equivalent to
\begin{equation*}
\vector{z}^{k+1} = \vector{0} \quad \text{ on } \Sigma,
\end{equation*}
for which no weak imposition by Lagrange multipliers is required.

\subsubsection{Offline stage}
\label{sec:ROM_2_offline}
During the offline stage, the solution of the high-fidelity problem $1_h$, $2.1_h$ and $2.2_h$ is first sought. Auxiliary unknowns $\vector{z}^{k+1}$ are then computed thanks to \eqref{eq:def_z}, for all $k = 0, ..., K - 1$. We then consider the following snapshot matrix
\begin{align*}
S_{\z} = [\underline{\mathbf{z}}^{1}_h | \hdots | \underline{\mathbf{z}}^{K}_h] \in \RR^{N_h^{\u} \times K},
\end{align*}
and, similarly to Section \ref{sec:ROM_1_offline}, retain the first $N_{\z}$ POD modes in the reduced space $V_N^{(2)}$. Reduced spaces for fluid pressure and structure displacement are defined as in Section \ref{sec:ROM_1_offline}, $Q_N^{(2)} := Q_N^{(1)}$ and $E_N^{(2)} := E_N^{(1)} = \text{span}(\{\phi_l\}_{l = 1}^{N^{\eta}})$. Moreover, harmonically extend each $\{\phi_l\}_{l = 1}^{N^{\eta}}$ to $\{\widehat{\phi}_l\}_{l = 1}^{N^{\eta}}$.

\subsubsection{Online stage}
\label{sec:ROM_2_online}
Similarly to Section \ref{sec:ROM_1_online}, a reduced-order approximation of the FSI problem is now obtained by means of a Galerkin projection over the reduced spaces $V_N^{(2)}, Q_N^{(2)}$ and $E_N^{(2)}$, respectively, that is: for any $k = 1, \hdots, K$
\begin{enumerate}
\item[1$_N^{(2)}$.] Explicit step (fluid viscous part), with change of variable for the fluid velocity: find ${\z}^{k+1}_N \in V_N^{(2)}$ such that:
\begin{align*}
\int_\Omega \frac{\rho_f}{\Delta t} \z^{k+1}_N \cdot \v_N \,d\x & + \int_\Omega 2 \mu_f \boldsymbol{\varepsilon}({\z}^{k+1}_N) : \grad\v_N \,d\x = \\
& \int_\Omega \frac{\rho_f}{\Delta t} \u^{k}_N \cdot \v_N \,d\x \, - \int_\Omega \grad p^k_N \cdot \v_N \,d\x, \\
& - \int_\Omega \frac{\rho_f}{\Delta t} D_t \widehat{\eta}^k_N\, \n \cdot \v_N \,d\x - \int_\Omega 2 \mu_f \boldsymbol{\varepsilon}(D_t \widehat{\eta}^k_N\, \n) : \grad\v_N \,d\x
\end{align*}
for all $\v_N \in V_N^{(2)}$. We note that boundary and \emph{interface} conditions on ${\z}^{k+1}_N$ are implicitly verified, since they are homogeneous and satisfied by any element in $V_N^{(2)}$. Finally, for the sake of the implicit step, we define $\u^{k+1}_N$ as $\vector{z}^{k+1}_N + D_t \widehat{\eta}^k_N\, \n$.
\item[2$_N^{(2)}$.] Implicit step: for any $j = 0, \hdots$, until convergence:
\begin{enumerate}
\item[2.1$_N^{(2)}$.] Fluid projection substep, with Robin boundary conditions: as in 2.1$_N^{(1)}$.
\item[2.2$_N^{(2)}$.] Structure substep: as in 2.1$_N^{(2)}$. Finally, at convergence, harmonically extend $\eta^{k+1}_N$ to $\widehat{\eta}^{k+1}_N$. Note that this can be easily carried out as the linear combination of the harmonically extended displacement basis $\{\widehat{\phi}_l\}_{l = 1}^{N^{\eta}}$. 
\end{enumerate}
\end{enumerate}

\begin{remark}[On efficient offline-online decoupling]
\label{rem:offline_online_decoupling_2}
Similarly to Remark \ref{rem:offline_online_decoupling_1}, an efficient offline-online decoupling can be obtained also in this case. Thanks to the definition of the harmonically extended displacements basis  $\{\widehat{\phi}_l\}_{l = 1}^{N^{\eta}}$ an efficient assembly of (e.g.) the right-hand side mass term $\int_\Omega \frac{\rho_f}{\Delta t} D_t \widehat{\eta}^k_N\, \n \cdot \v_N \,d\x$ can be obtained. We remark that this accounts for a negligible additional offline cost (solution of \reviewerA{$N^\eta$} harmonic extension problems) and no additional online cost, since the extension of $\eta^{k+1}_N$ to $\widehat{\eta}^{k+1}_N$ does not actually require the solution of reduced problem, but rather a linear combination of $\{\widehat{\phi}_l\}_{l = 1}^{N^{\eta}}$ \reviewerA{once the coefficients of the structural unknown have been computed}.

\end{remark}

\section{Numerical comparison}
\label{sec:numerical_comparison}

In this section we summarize the numerical results obtained by the  proposed reduced-order models. The values of constitutive and geometrical parameters are summarized in Table \ref{tab:constpar} \cite{Formaggia2001}. The domain has been discretized with a $120 \times 10$ structured mesh, while the time-step is $\Delta t = 10^{-4} \text{ s}$. The final time is $T = 0.13 \text{ s}$, so that $K = 1300$. $T$ has been chosen to simulate the pressure wave propagation, just before wave reflection occurs. Numerical simulations are carried out using \emph{RBniCS} \cite{HesthavenRozzaStamm2015,RBniCS}, an open-source reduced order modelling library developed at SISSA mathLab, built on top of \emph{FEniCS} \cite{LoggMardalEtAl2012a}.

\begin{table}
\centering
\begin{tabular}{|l|l||l|l||l|l|}\hline
$\rho_f$ & $1 \text{ g/cm}^3$ & $\mu_f$ & $0.035 \text{ Poise}$ & $\rho_s$ & $1.1 \text{ g/cm}^3$ \\\hline
 $E_s$ & $0.75 \times 10^6 \text{ dyn/cm}^2$ & $\nu_s$ & $0.5$ & $c_1$ & $\frac{h_s E_s}{h_f^2 \ (1 - \nu_s^2)}$ \\\hline
$c_0$ & $\frac{h_s E_s}{2\ (1+\nu_s)}$ & $L$ & $6 \text{ cm}$ & $h_f$ & $0.5 \text{ cm}$ \\\hline 
$h_s$ & $0.1 \text{ cm}$ & $p_{in}$ & $10^4 (1-\cos(\frac{2 \pi t}{0.005})) \boldsymbol{1}_{t < 0.005} \text{ dyn/cm}^2$ & 
$p_{out}$ & $0 \text{ dyn/cm}^2$\\\hline
\end{tabular}
\caption{Constitutive parameters for test case (from \cite{Formaggia2001}).}
\label{tab:constpar}
\end{table}

\begin{figure}
\centering
\subfloat[POD singular values.]{\includegraphics[width=0.48\textwidth]{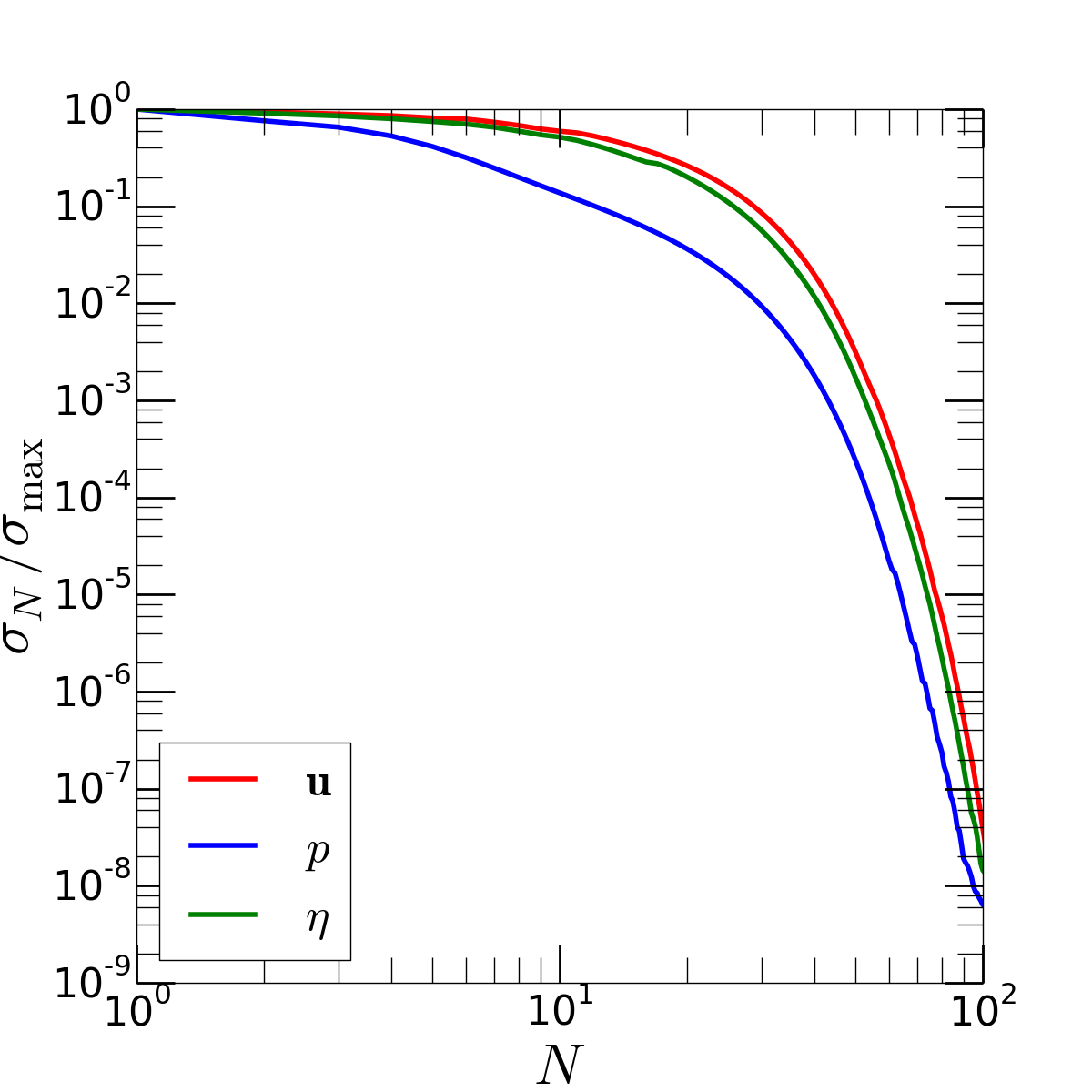}}\;
\subfloat[POD retained energy.]{\includegraphics[width=0.48\textwidth]{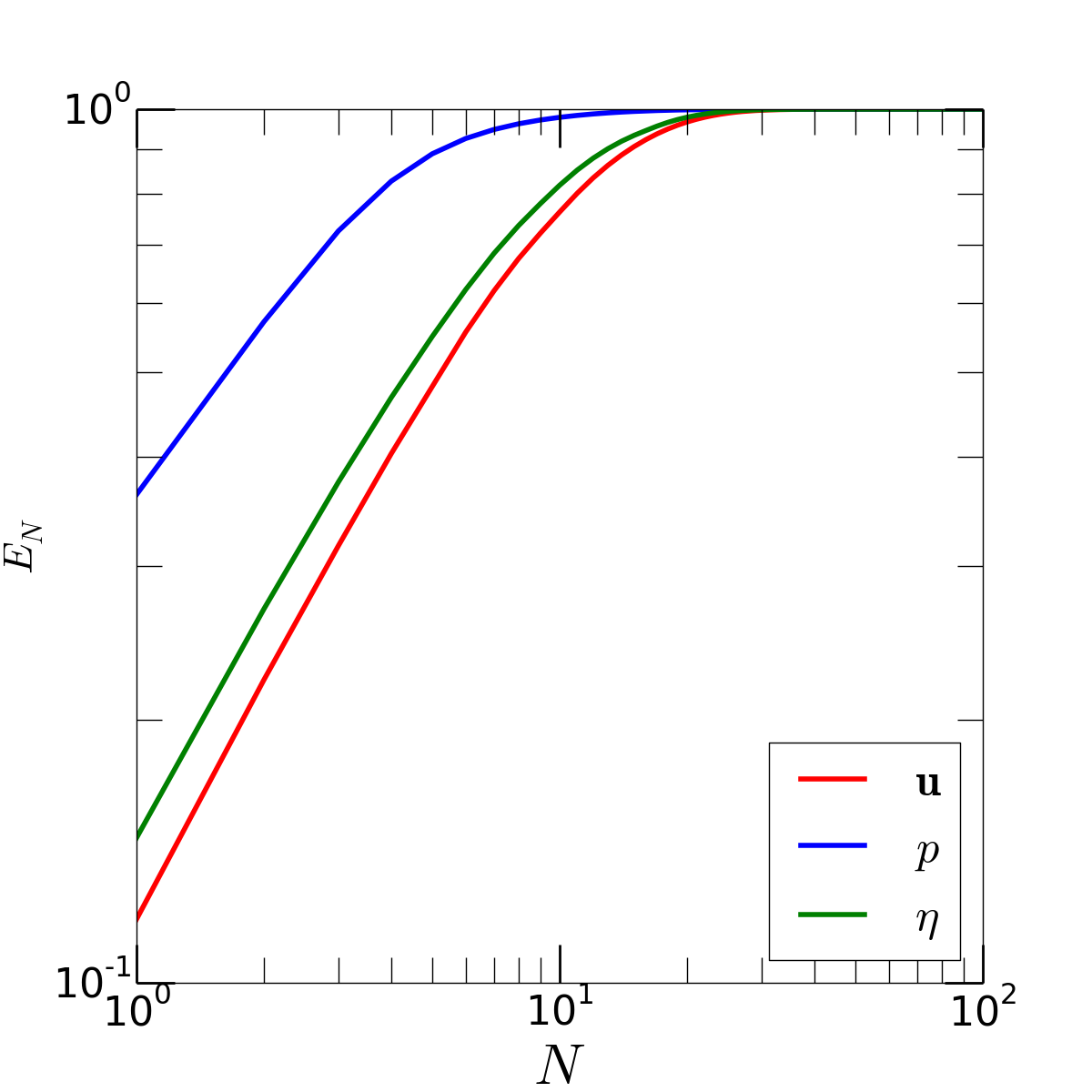}}
\caption{Results of the offline stage of FSI ROM 1: POD singular values and retained energy as a function of the number $N$ of POD modes for fluid velocity, fluid pressure, and solid displacement.}
\label{fig:POD_ROM_1}
\end{figure}

\begin{figure}
\centering
\subfloat[Error analysis.]{\includegraphics[width=0.48\textwidth]{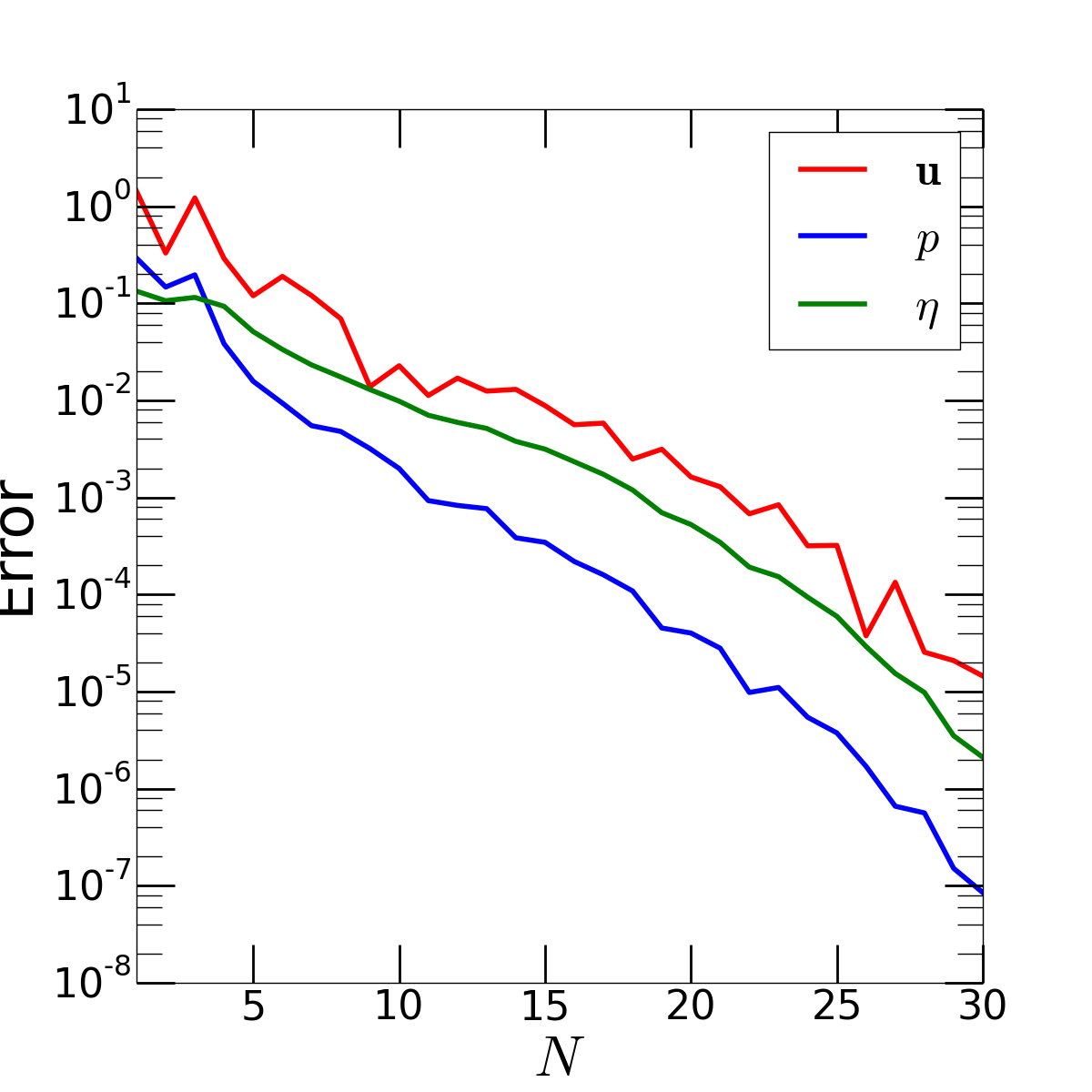}\label{fig:error_ROM_1}}\;
\subfloat[Speedup analysis.]{\includegraphics[width=0.48\textwidth]{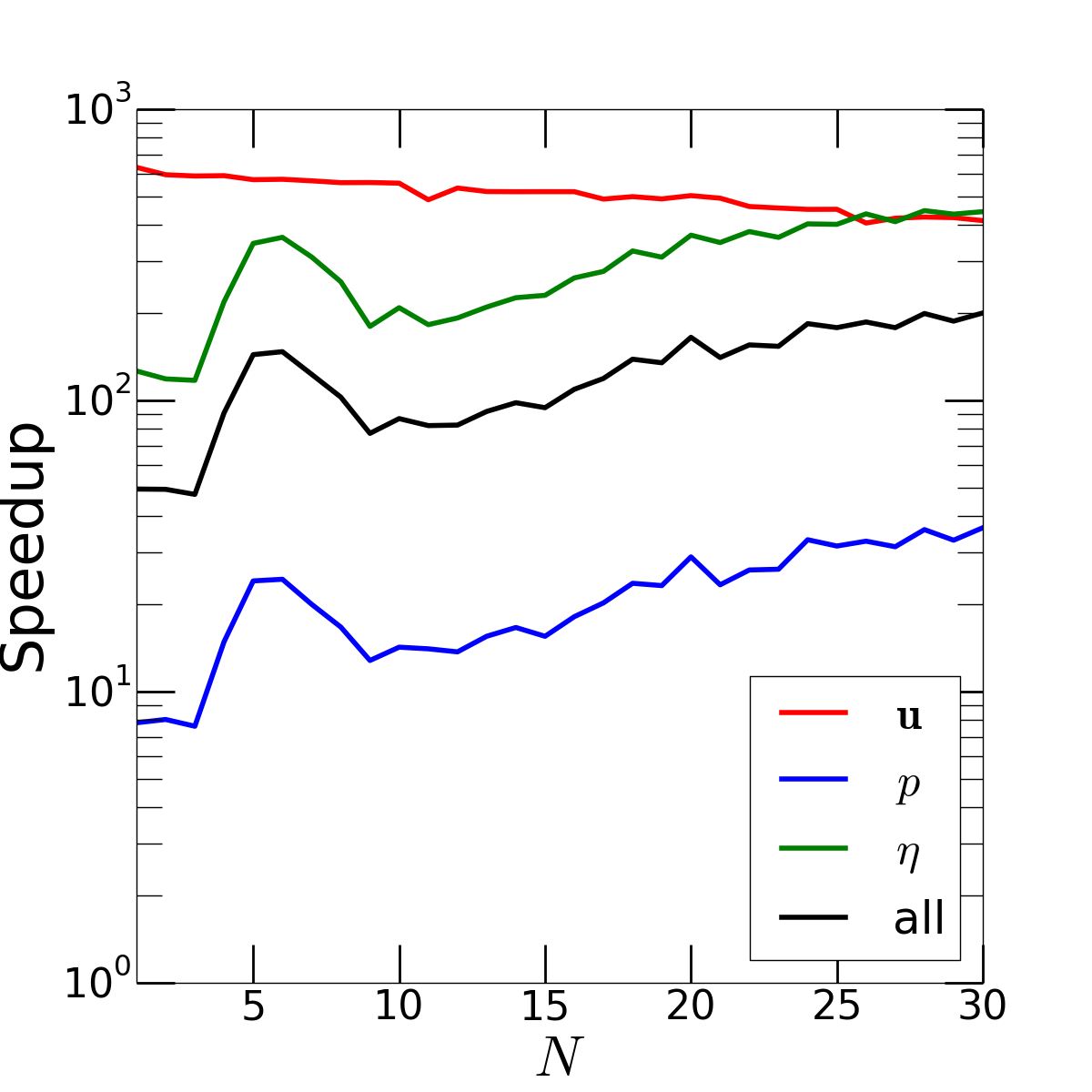}\label{fig:speedup_ROM_1}}
\caption{Error analysis and speedup analysis of FSI ROM 1, as a function of the number $N$ of POD modes for fluid velocity, fluid pressure, and solid displacement.}
\end{figure}

\begin{figure}
\centering
\subfloat[Speedup analysis of FSI ROM 2.]{\includegraphics[width=0.48\textwidth]{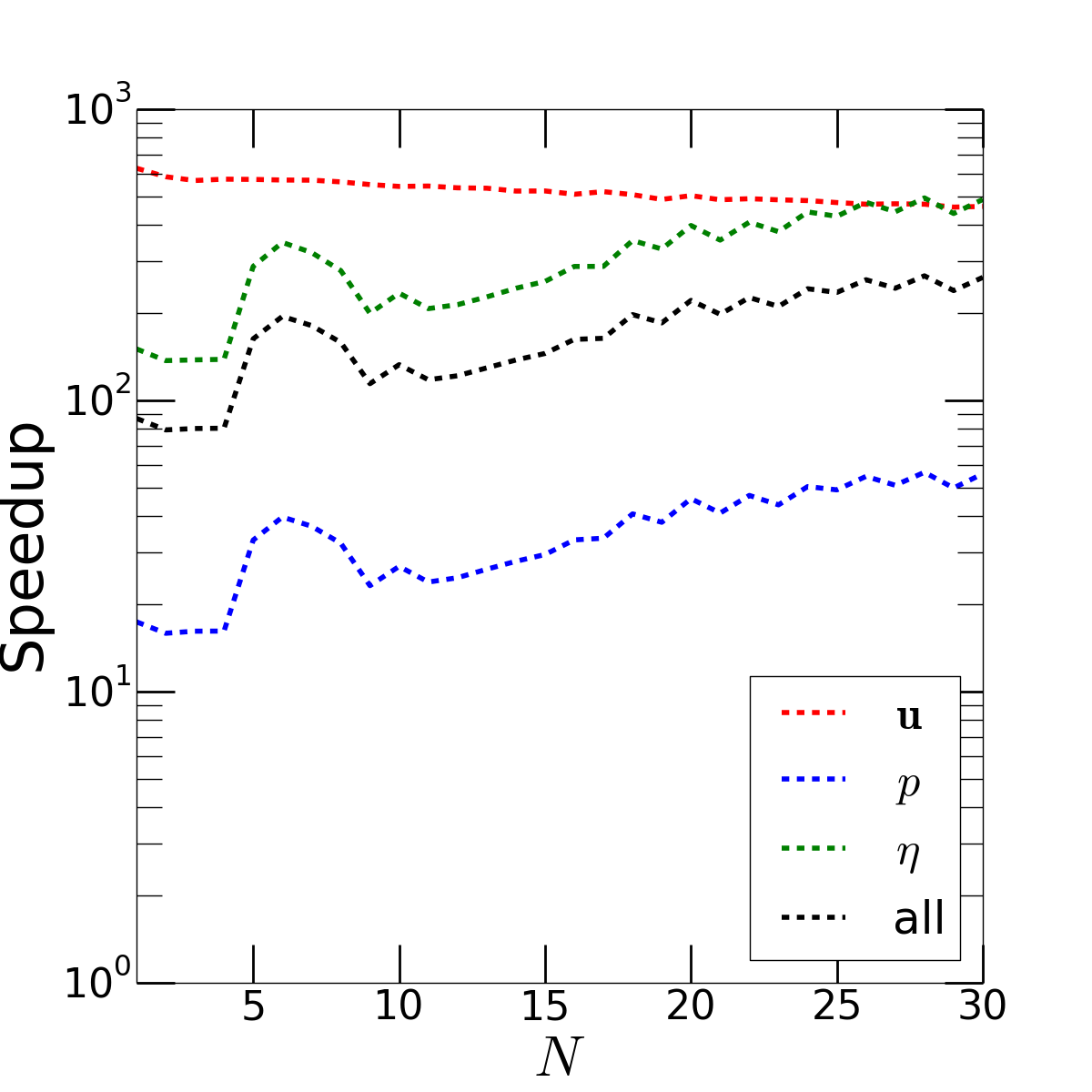}}\;
\subfloat[Comparison of the speedup for fluid velocity unknown of FSI ROMs 1 and 2.]{\includegraphics[width=0.48\textwidth]{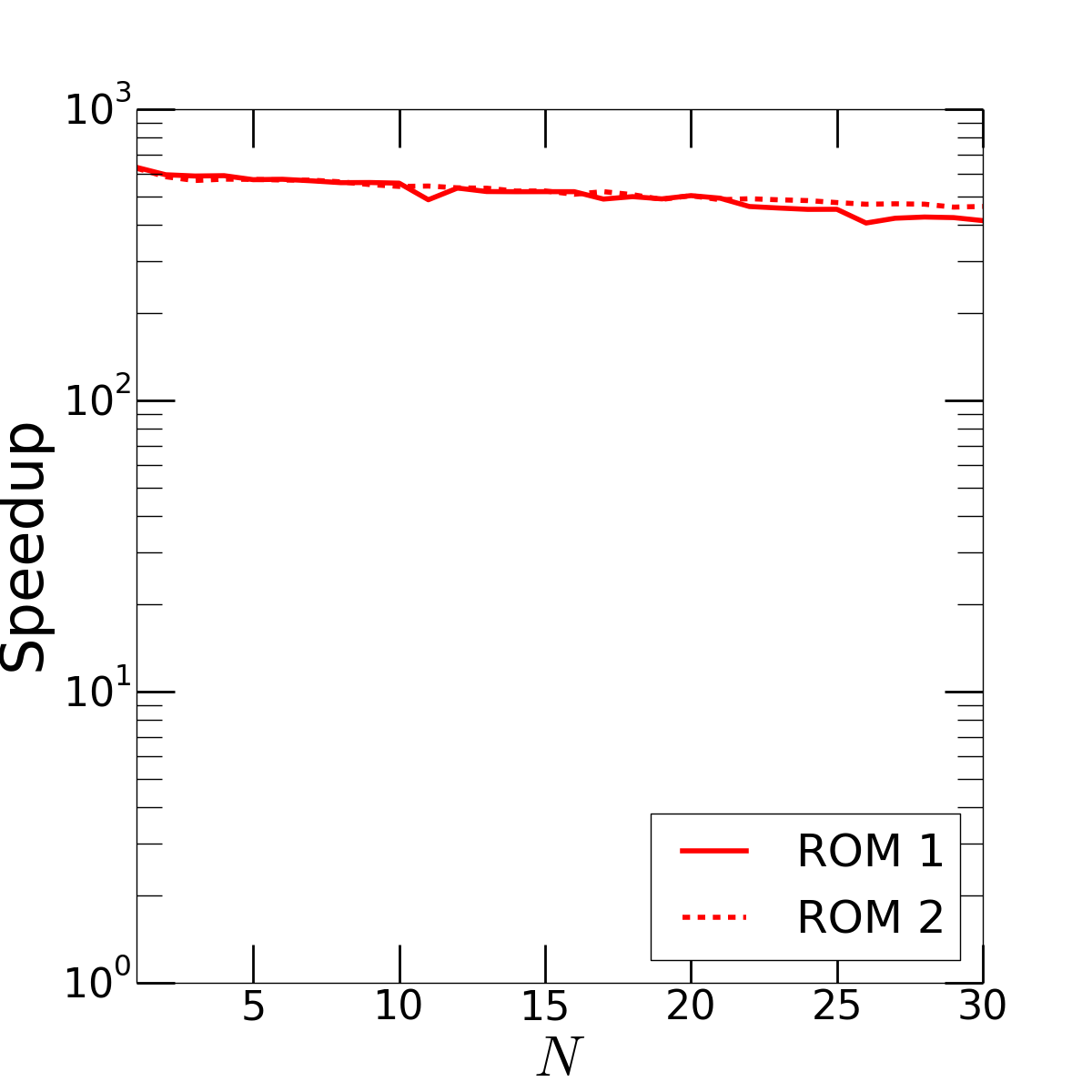}}
\caption{Speedup analysis of FSI ROM 2, as a function of the number $N$ of POD modes for fluid velocity, fluid pressure, and solid displacement.}
\label{fig:speedup_ROM_2}
\end{figure}

\vspace{\baselineskip}
\noindent \textbf{FSI ROM 1:} Fig. \ref{fig:POD_ROM_1} shows the POD singular values and retained energy as a function of the number $N$ of POD modes for FSI ROM 1. It can be noticed that the decay of the singular values of fluid velocity and structure displacement is slower than the one of fluid pressure; moreover, the first POD mode of fluid pressure retains a larger energy ($\approx 36\%$) than the first modes of structure displacement ($\approx 15\%$) and fluid velocity ($\approx 12\%$). Accordingly, the (relative) error analysis (Fig. \ref{fig:error_ROM_1}) shows that the reduced solution converges to the high-fidelity one, and that (except for small $N$), the relative error on the pressure is smaller than the displacement one, which is smaller than the velocity relative error. Employing only $N=30$ POD modes out of the $K=1300$ snapshots, velocity, displacement and pressure relative errors are of the order of $10^{-4}$, $10^{-5}$ and $10^{-7}$, respectively. Figure \ref{fig:speedup_ROM_1} shows that the overall speedup for large values of $N$ is of at least two order of magnitudes. In particular, speedup for the \reviewerA{explicit} step is approximately constant for increasing $N$, \reviewerA{while the speedup for the implicit steps is increasing for larger $N$ since the number of required iterations for the implicit step decreases with $N$ (see Fig. \ref{fig:speedup_ROM_2}b and c for the maximum and average number of iterations, respectively).}

\begin{figure}
\centering
\subfloat[POD singular values.]{\includegraphics[width=0.48\textwidth]{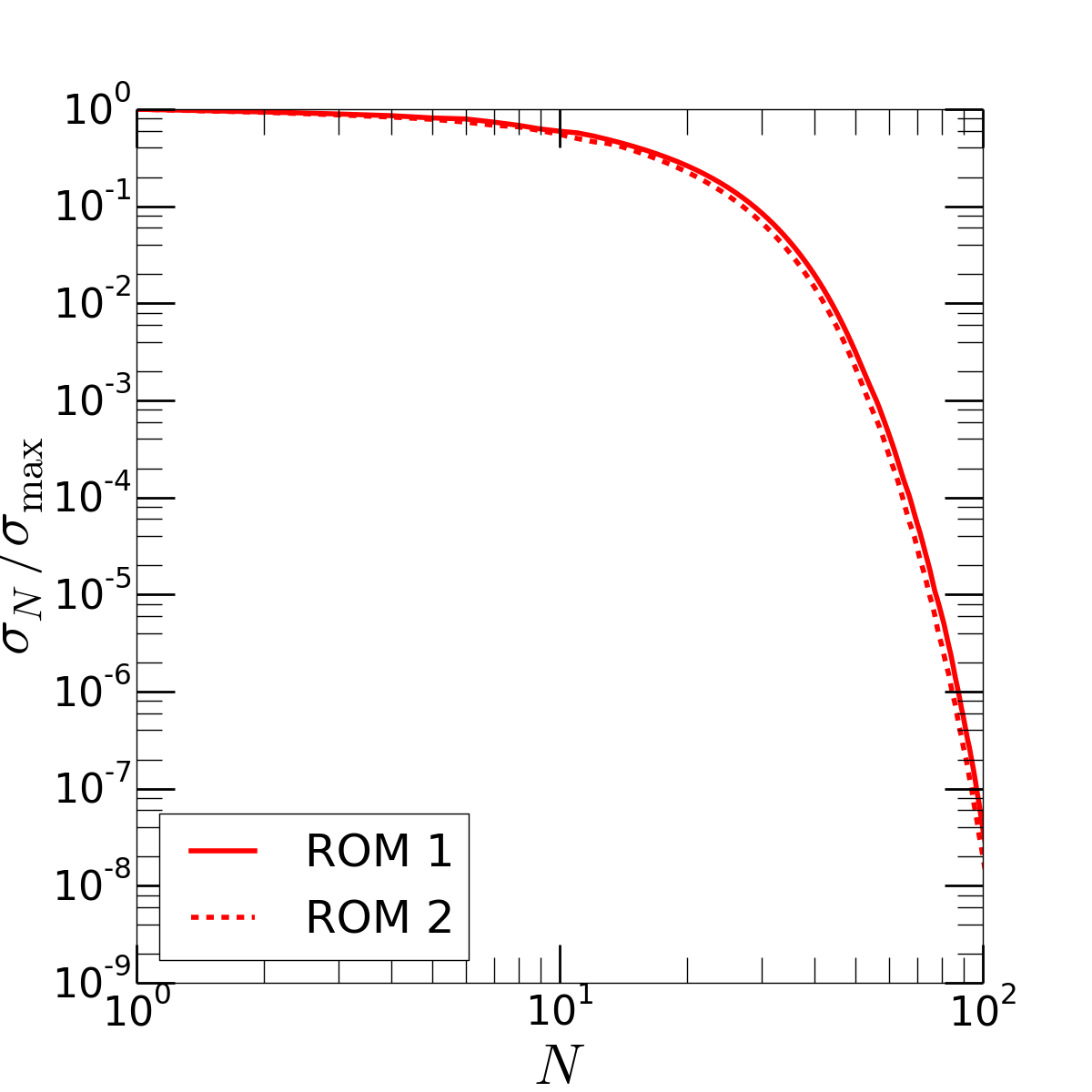}}\;
\subfloat[POD retained energy.]{\includegraphics[width=0.48\textwidth]{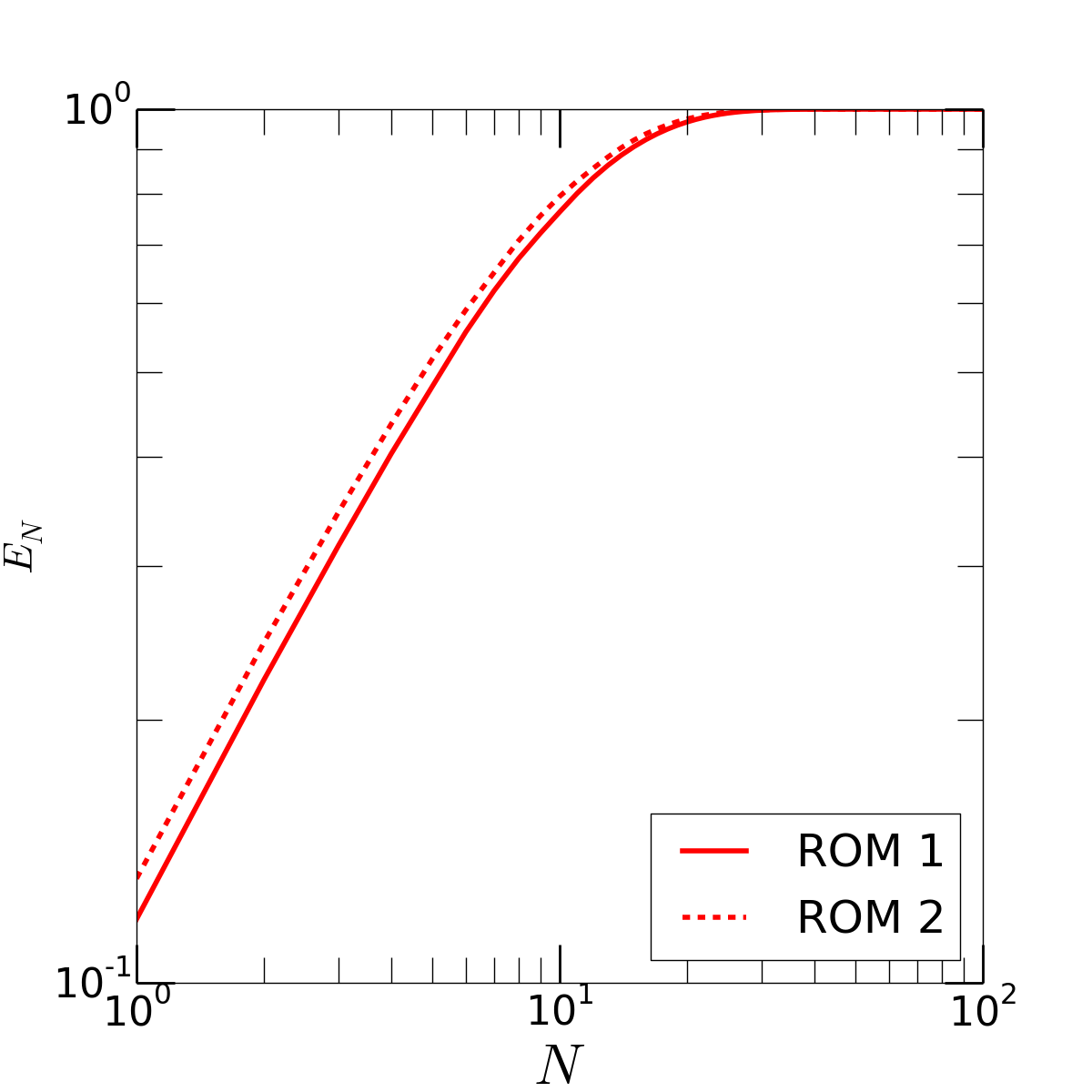}}
\caption{Comparison of the offline stage of FSI ROMs 1 and 2: POD singular values and retained energy as a function of the number $N$ of POD modes for fluid velocity $\u$ (FSI ROM 1) and auxiliary fluid velocity $\vector{z}$ (FSI ROM 2).}
\label{fig:POD_ROM_2}
\end{figure}

\begin{figure}
\centering
\subfloat[Condition number of explicit step.]{\includegraphics[width=0.48\textwidth]{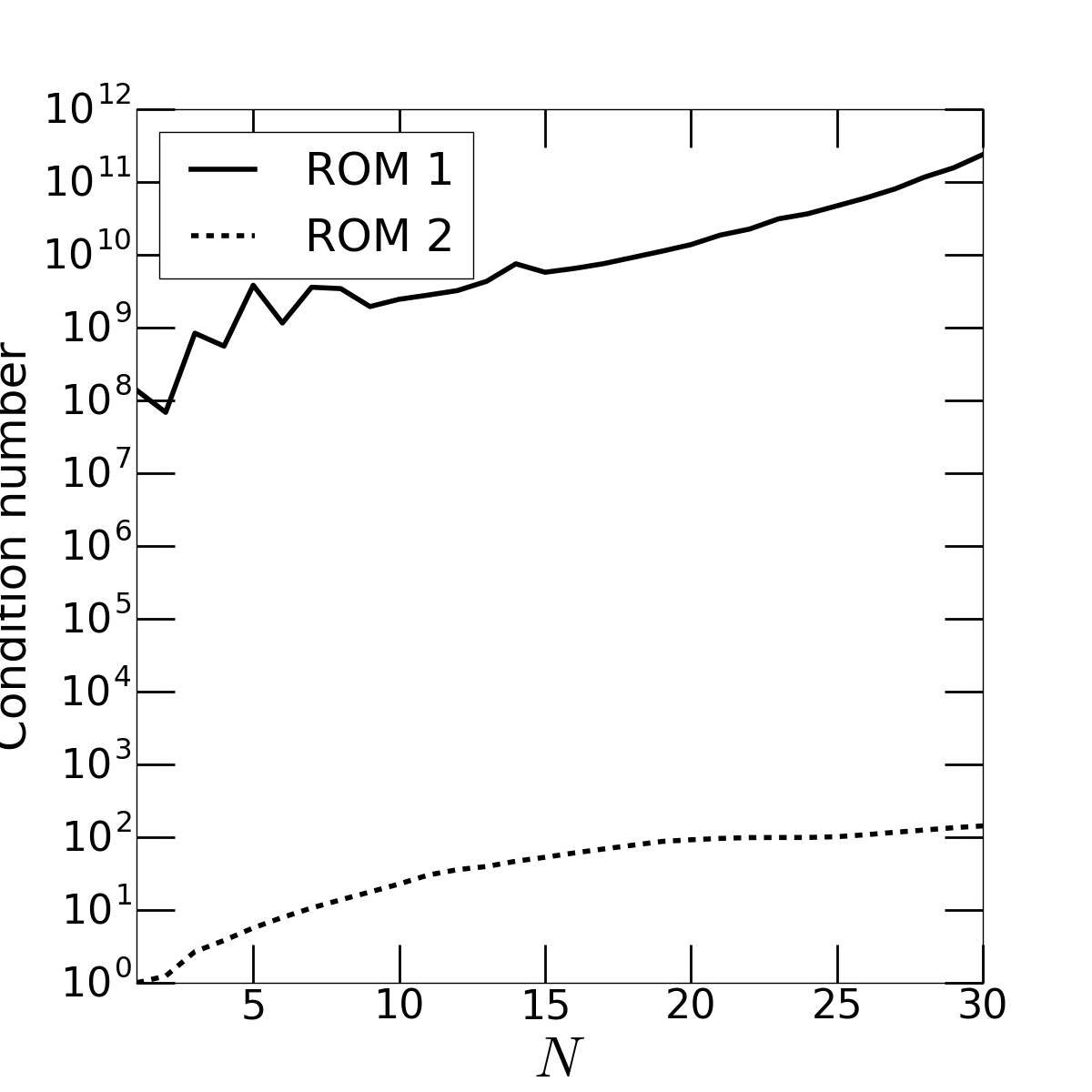}\label{fig:cond1_ROM_2}}\\
\subfloat[Maximum number of iterations of implicit step.]{\includegraphics[width=0.48\textwidth]{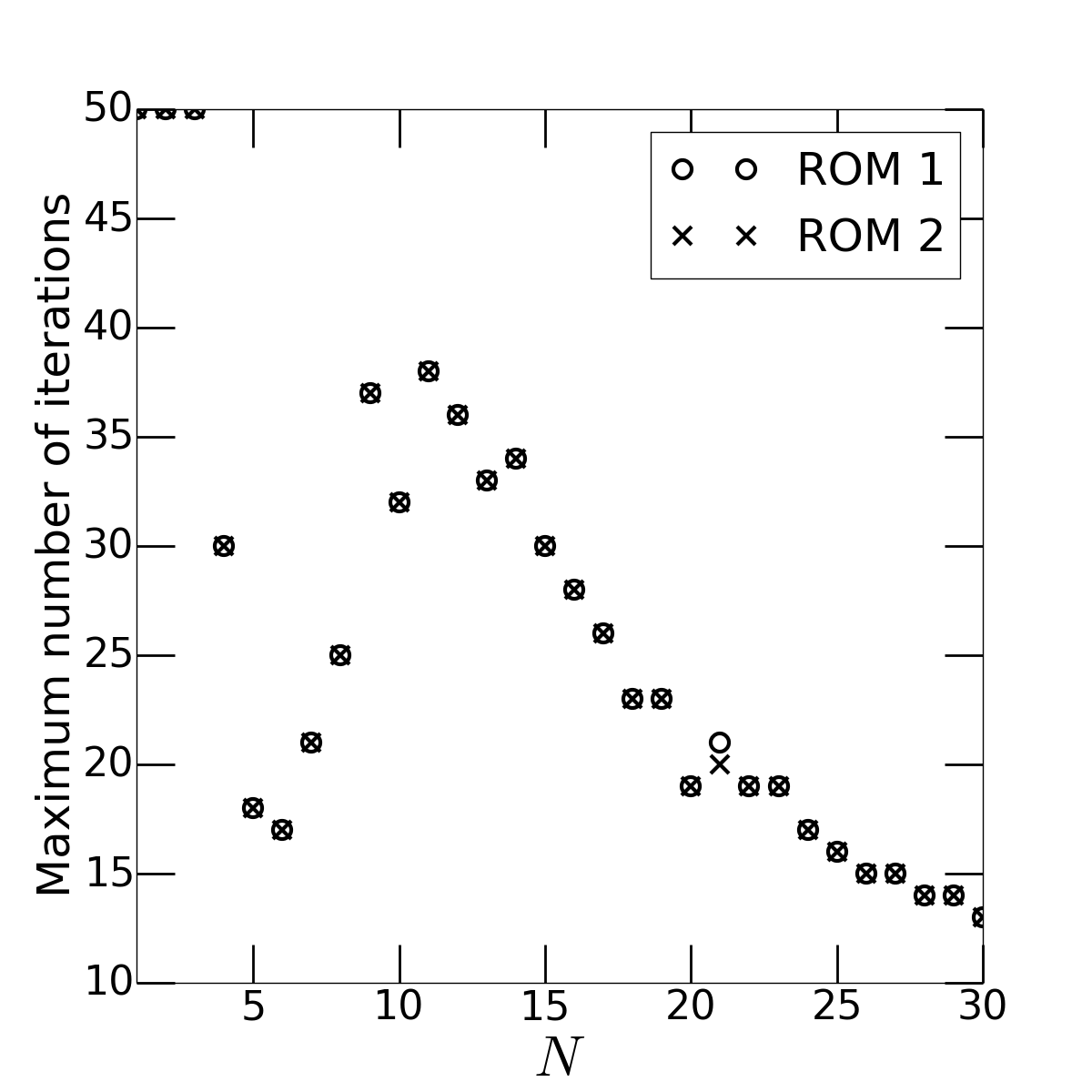}\label{fig:itmax_ROM_2}}\;
\subfloat[\reviewerA{Average number of iterations of implicit step.}]{\includegraphics[width=0.48\textwidth]{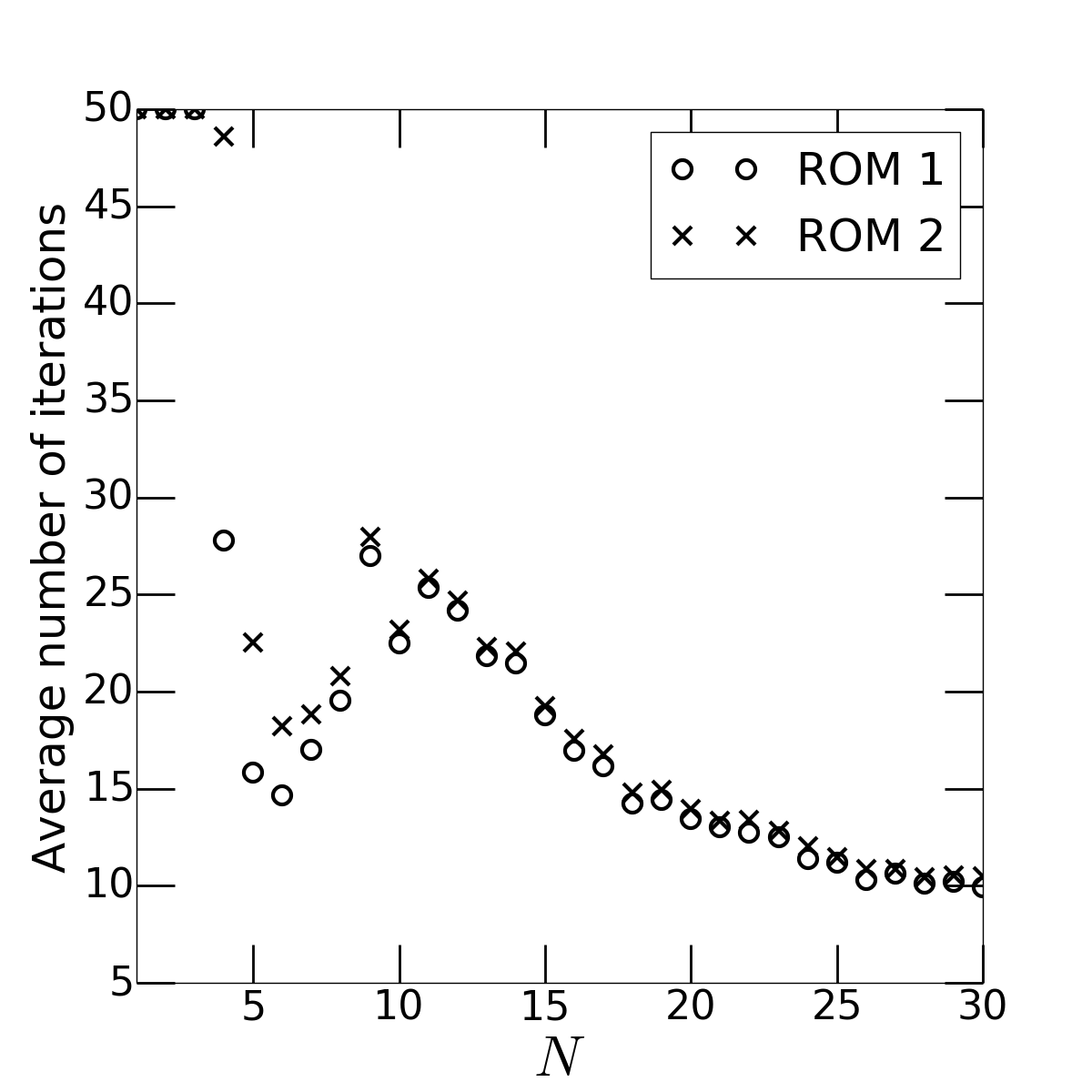}\label{fig:itave_ROM_2}}
\caption{Comparison of the condition number of the left-hand side matrix of $1_N^{(1)}$ and $1_N^{(2)}$ (top) and \reviewerA{of the maximum (bottom, left) and average (bottom, right) number of iterations required by $2_N^{(1)}$ and $2_N^{(2)}$}, as a function of the number $N$ of POD modes for fluid velocity, fluid pressure, and solid displacement.}
\end{figure}

\begin{figure}
\centering
\subfloat[Error analysis of FSI ROM 2.]{\includegraphics[width=0.48\textwidth]{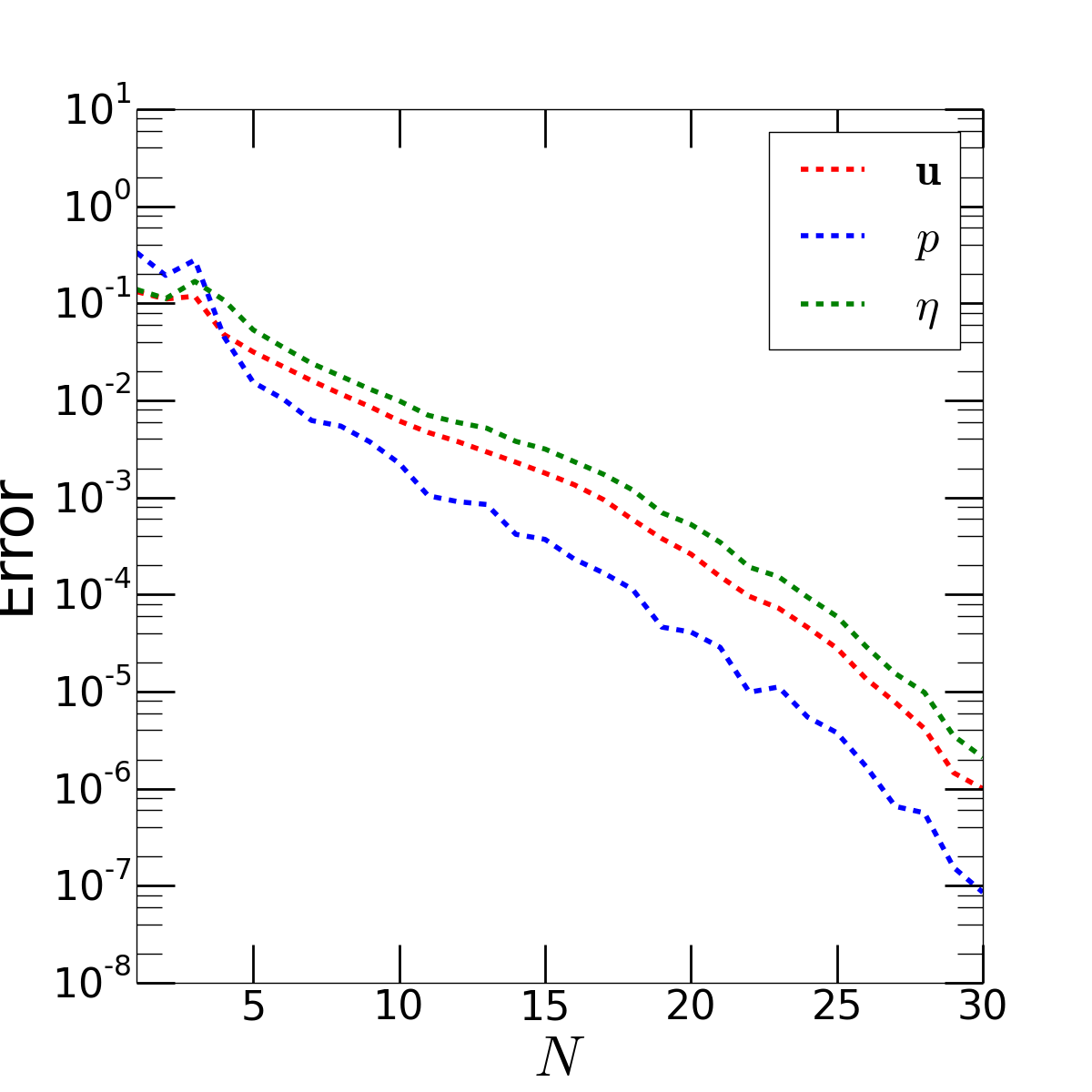}}\;
\subfloat[Comparison of the error for fluid velocity unknown of FSI ROMs 1 and 2.]{\includegraphics[width=0.48\textwidth]{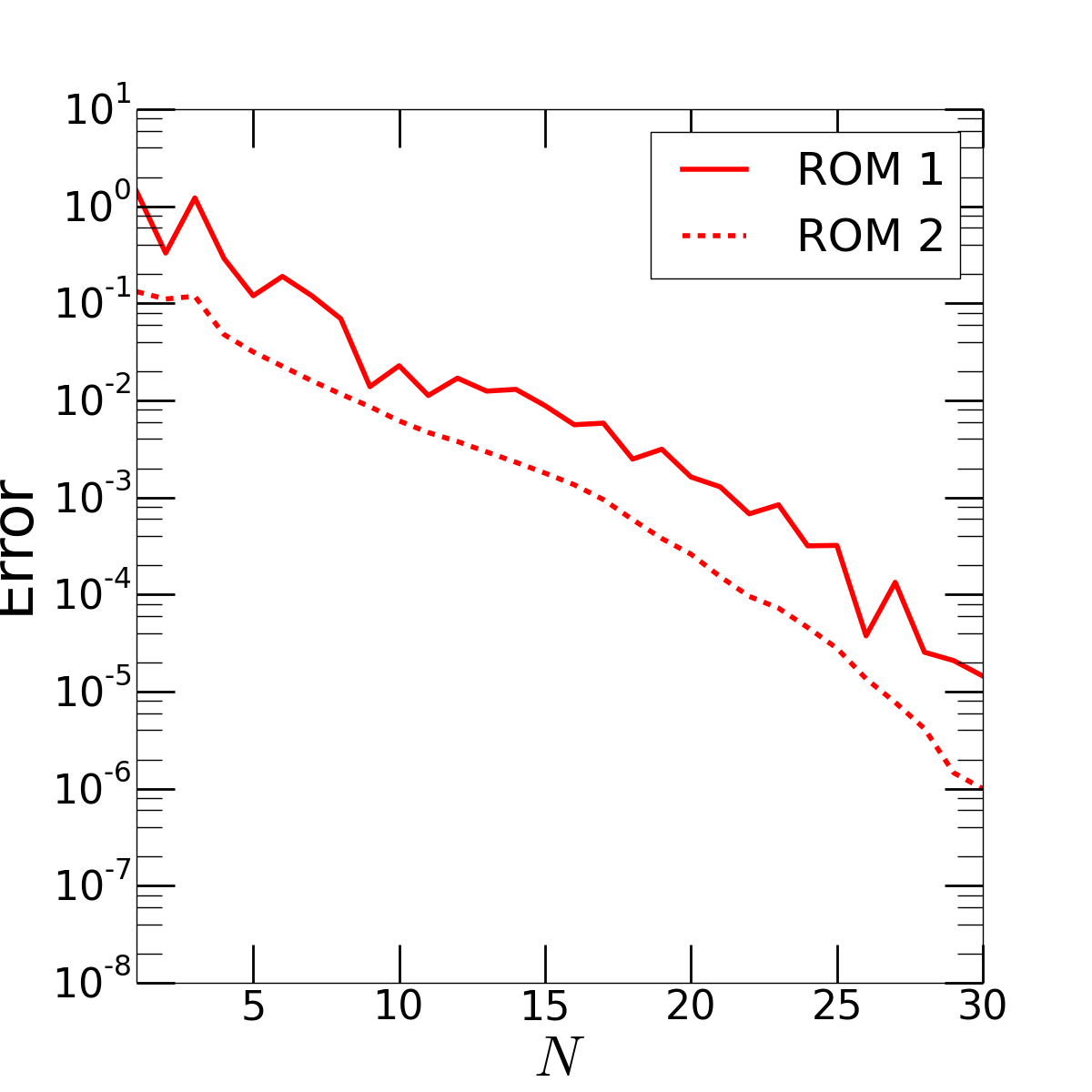}}
\caption{Error analysis of FSI ROM 2, as a function of the number $N$ of POD modes for fluid velocity, fluid pressure, and solid displacement.}
\label{fig:error_ROM_2}
\end{figure}

\begin{figure}
\centering
\includegraphics[width=0.48\textwidth]{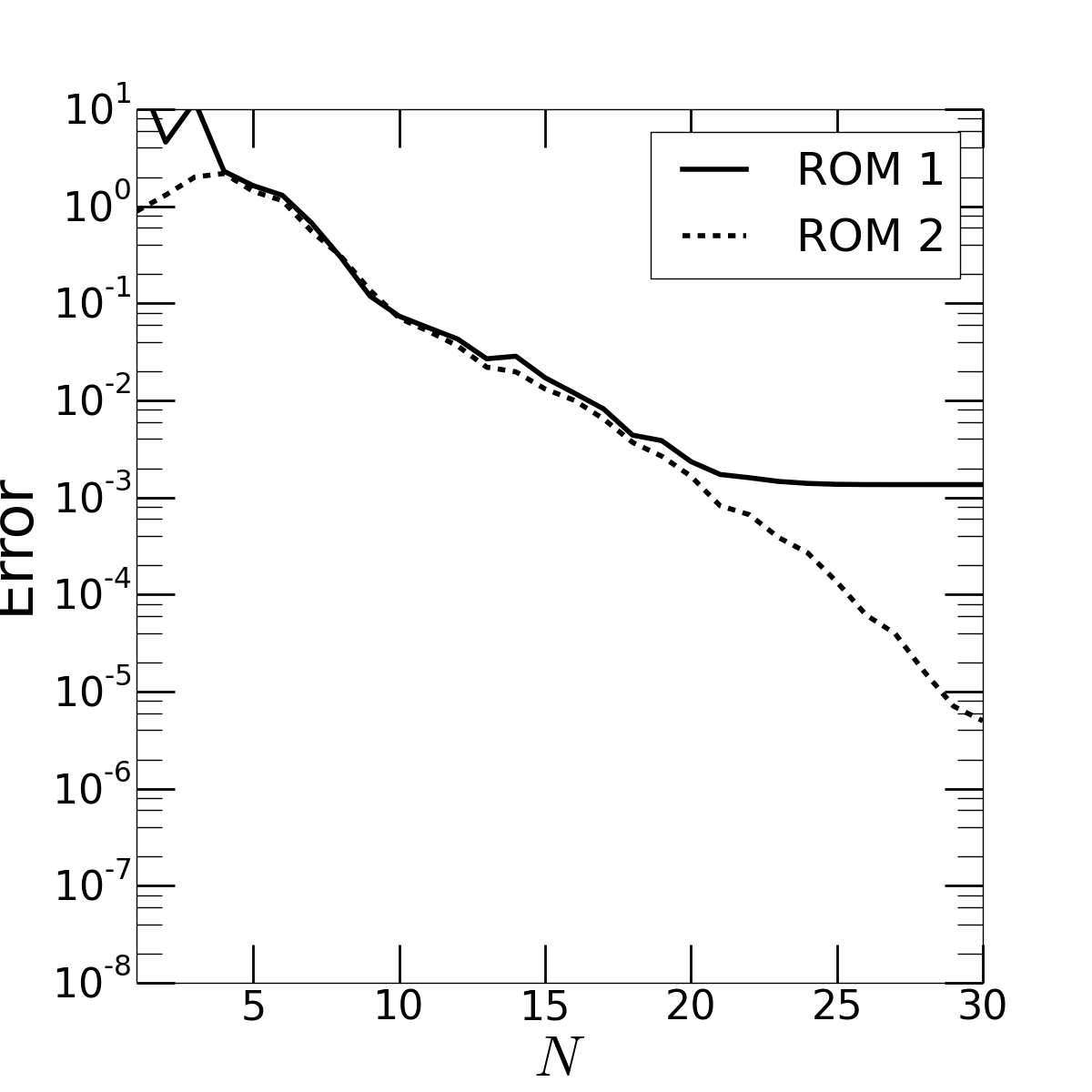}
\caption{\reviewerA{Error analysis of the interface stress for FSI ROMs 1 and 2, as a function of the number $N$ of POD modes for fluid velocity, fluid pressure, and solid displacement.}}
\label{fig:error_ROM_stress}
\end{figure}

\vspace{\baselineskip}
\noindent \textbf{FSI ROM 2:} Fig. \ref{fig:POD_ROM_2} shows a comparison of the offline stage for ROMs 1 and 2. From Fig. \ref{fig:POD_ROM_2}b it can be noticed that, after the change of variable, the first auxiliary velocity $\z$ POD mode of FSI ROM 2 retains approximately $1.5\%$ more energy than the first velocity $\u$ POD mode. The remaining variables are omitted since the pressure and displacement bases are the same as in FSI ROM 1. Moreover, Fig. \ref{fig:cond1_ROM_2} shows that the left-hand side matrix of the fluid viscous step $1_N^{(2)}$ (FSI ROM 2) is characterized by a condition number which is, for all $N$, at least 10 order of magnitude smaller than the one of $1_N^{(1)}$ (FSI ROM 1). The combination of these two remarks justifies the improvement in the error analysis for the velocity variables, shown in Fig. \ref{fig:error_ROM_2}b. On average, the relative error on the velocity unknown obtained by FSI ROM 1 is seven times the one obtained by FSI ROM 2. Relative errors for the remaining unknowns are omitted because they are comparable among FSI ROMs 1 and 2. \reviewerA{Moreover, we show in Fig. \ref{fig:error_ROM_stress} the error analysis for the interface stress. The plot clearly shows that FSI ROM 2 provides a better approximation of the interface stress for $N > 20$.}
Online performance (Fig. \ref{fig:speedup_ROM_2}) are comparable to the ones obtained by FSI ROM 1. \reviewerA{This is due to the fact that (i) the number of iterations to reach convergence in step $2_N^{(1)}$ and $2_N^{(2)}$ are the same as maximum values (Fig. \ref{fig:itmax_ROM_2}) and comparable on average (Fig. \ref{fig:itave_ROM_2})}, and (ii) the time to solve the explicit step does not depend strongly on $N$ (Fig. \ref{fig:speedup_ROM_2}b), even though step $1_N^{(2)}$ (FSI ROM 2) requires the solution of a linear system of size $N \times N$ (rather than $2 N \times 2 N$ for step $1_N^{(1)}$ (FSI ROM 1)).

\section{Conclusions and perspectives}
\label{sec:conclusion}
Two semi-implicit reduced-order models for FSI problems have been proposed in this work, based on a POD--Galerkin approximation of an operator splitting semi-implicit high-fidelity scheme. FSI ROM 1 is a standard Galerkin projection over the reduced spaces generated by POD. No supremizer enrichment is required, thanks to the operator splitting approach. Even though FSI ROM 1 shows good performance in terms of error analysis, its major drawbacks (when compared to FSI ROM 2) are related to the weak imposition of \eqref{eq:bc_strongly} by Lagrange multipliers. Numerical results of the previous section have shown that this is detrimental for several aspects of the ROM: increased system dimension of the fluid explicit step, increased condition number of the fluid explicit step, increased error for the velocity. FSI ROM 2, instead, stems from the idea that \eqref{eq:bc_strongly} can be easily imposed in a reduced-order framework by performing the change of variable \eqref{eq:def_z}. In this way, all the detrimental effects of FSI ROM 1 are remedied. Moreover, better properties in terms of POD retained energy are also obtained. The combination of these two factors results in a better approximation of the fluid velocity. In particular, in FSI ROM 2, we try to separate in the fluid velocity the fluid-structure interaction component from the pure fluid part. 
\reviewerA{Even though FSI ROM 1 suffers several drawbacks when compared to FSI ROM 2 (especially for what concerns the increased condition number), the approach proposed by FSI ROM 1 can be more easily integrated with existing reduced order modelling capabilities for fluid problems, since it does not require to change the existing computations of fluid velocity basis functions. Nevertheless, especially for FSI ROM 1, a more detailed analysis of enrichment procedures shall be carried out to further investigate the stability of the resulting reduced problem, due to two saddle point structures to be taken into account (see Remarks 3 and 4).}
Future work will concern ROMs that are better able to face the hyperbolic nature of the problem. 
A more efficient separation at the reduced-order level of the parabolic and hyperbolic components of the system may decrease the number of basis functions required to obtain an accurate reduced description of the FSI problem.


\section*{Acknowledgements}
We acknowledge the support by European Union Funding for Research and Innovation -- Horizon 2020 Program -- in the framework of European Research Council Executive Agency: H2020 ERC CoG 2015 AROMA-CFD project 681447 ``Advanced Reduced Order Methods with Applications in Computational Fluid Dynamics'' and the PRIN project ``Mathematical and numerical modelling of the cardiovascular system, and their clinical applications''. We also acknowledge the INDAM-GNCS project ``Tecniche di riduzione della complessit\`{a} computazionale per le scienze applicate'' and INDAM-GNCS young researchers project ``Numerical methods for model order reduction of PDEs''.

{
\small
\bibliography{bib/library} 
\bibliographystyle{spmpsci}
}

\end{document}